\documentclass{svjour2}                    
\smartqed  
\usepackage{graphicx}
\usepackage{latexsym, amssymb, amsmath}
\newcommand{\lf}{\left}
\newcommand{\rt}{\right}

\newcommand{\be}{\begin{equation}}
\newcommand{\ee}{\end{equation}}

\newcommand{\al}{\alpha}
\newcommand{\ep}{\epsilon}
\newcommand{\de}{\delta}

\newcommand{\la}{\lambda}

\newcommand{\tht}{\theta}
\newcommand{\half}{\frac{1}{2}}
\newcommand{\tP}{\tilde{P}}
\newcommand{\tQ}{\tilde{Q}}
\newcommand{\moP}{\tP}
\newcommand{\CL}{\mathcal{L}}
\DeclareMathOperator{\mE}{E}
\DeclareMathOperator{\mP}{P}
\DeclareMathOperator{\mL}{\CL}
\newcommand{\p}{\partial}
\newcommand{\bF}{\bar{f}} 

\newcommand{\tbf}{\textbf}
\newcommand{\n}{\noindent}

\begin{document}

\title{Reinforced Brownian Motion on the Half-Line }
\subtitle{}

\titlerunning{}        

\author{Jerome K. Percus \and
        Ora E. Percus
	}


\institute{Jerome K. Percus \at
              Courant Institute/NYU, 251 Mercer Street, New York, NY 10012 \\
              Tel.: (212) 998-3130\\
              Fax: (212) 995-4120 \\
              \email{percus@cims.nyu.edu}           
\and
Ora E. Percus \at
Courant Institute/NYU, 251 Mercer Street, New York, NY 10012
}

\date{Received: date / Accepted: date}

\maketitle

\begin{abstract}
We analyze the Brownian Motion limit of a prototypical unit step reinforced 
random-walk on the half line.  A reinforced random walk is one which changes 
the weight of any edge (or vertex) visited to increase the frequency of return 
visits.  The generating function for the discrete case is first derived for the
joint probability distribution of $S_N$ (the location of the walker at the 
$N^{th}$ step) and $A_N$ the maximum location the walker achieved in $N$ steps.  
Then the bulk of the analysis concerns the statistics of the limiting Brownian 
walker, and of its ``environment'', both parametrized by the amplitude
of the reinforcement.

\keywords{walk on half-line, reinforced random walk, Brownian Motion limit}
\end{abstract}

\section{Introduction}

A random walk on a lattice is termed edge (or vertex) reinforced\cite{bib1}
if it changes the parameters of any edge [or vertex] visited to increase the
frequency of return visits.  Thus, in principle, the state space of the walk
is that of the walker and its environment.  ``Long time potentiation'' of
neural signals through a network of synaptic junctions serves perhaps as a 
realistic conceptual model.  For even a qualitative appreciation of what is
involved, one must be much more specific.  In this paper, we will greatly 
reduce our scope, without losing sight of the distant target.  
To start, we restrict our initial attention to a Markovian walk on a 
one-dimensional integer lattice (reflected at the origin) with only 
next neighbor steps $X_i = \pm 1$.  Indexing by step number $N$, the
state of the walker is then specified by the location after $N$ steps
$S_N=\sum^N_{i=1} X_i$ and we also choose to define the state of lattice
by the position $A_N$, the highest value attained by the walker during 
the first $N$ steps.
Pictorially, one could imagine that any edge traversed by the walker
was permanently marked in bold print.  The walk on the joint space, 
$(S_N, A_N)$ is Markovian, and so it was convenient to investigate its
properties and then those of its two marginals.  Including reflection
at the origin, the probability distribution $P_N(x,a) = P\{ S_N=x, A_N=a\}$
was readily seen to satisfy
\begin{align*}
P_1 (x,a) &= \de_{x_1\, 1}\, \de_{a,1} \\
P_{N+1} (x,a) &= \half \lf( 1+ \de_{x_1\,1} - \de_{x, a+1} \rt) 
P_N \lf(x-1, a\rt) \\
&+ \half \: P_N \lf(x+1,a\rt) + \half \:P_N \lf(a-1, a-1\rt) 
\de_{x,a} \lf(1-\de_{a_1}\rt)
\end{align*}
Standard analysis led to an explicit expression for the 2 variable
generating function (g.f)
$$
\mP^{\approx} \lf(\la, u,a\rt)  = \sum^\infty_{N=1} \sum^\infty_{x=0}
P_N (x,a) \la^N\, u^x.
$$
From this the 2 marginal $g,fs$ are easily obtained.
Note that for recovering the marginal distribution function for $S_N$
it is not necessary to use the $g,f$ since it can be easily obtained,
see \cite{bib2}  eq 2.7 by familiar and routine argument.  Obtaining
$P\{A_N=a\}$ is neither familiar nor routine and is not Markovian as
well.  This was found as a rather complicated infinite series.  However
it was observed that the asymptotic limit of  $aP\{ A_N=a\}$ as 
$N\to\infty$ at fixed $\gamma$ where $\frac{a^2}{N} = \gamma \frac{\pi}{2}
+O\lf(\frac{1}{a}\rt)$ could be simplified to a series which converges 
very rapidly for $\gamma \ge 1$ and we had the asymptotic result
\begin{gather*}
(1-\al) \,\sqrt{\gamma}\, e^{-\frac{\pi \gamma}{4}} \le \lim_{N\to\infty}
aP\{ A_N=a\} \le \sqrt{\gamma}\, e^{-\frac{\pi\gamma}{4}}
\\
\text{for}\quad\gamma\ge 1 \qquad \text{where}\qquad \al=3e^{-2\pi\gamma}\le
0.0056.
\end{gather*}
The above series could also be Poisson transformed to another series which 
converges very rapidly for $\gamma \le 1$, yielding the asymptotic result
\begin{gather*}
\lf(1-\al'\rt) \frac{1}{\gamma}\:e^{-\frac{\pi}{4\gamma}} \le 
\lim_{N\to\infty} \; aP\{ A_N=a\} \le \frac{1}{\gamma}\:
e^{-\frac{\pi}{4\gamma}} \\
\text{for}\quad \gamma \le 1 \qquad \text{ where } \qquad \al'=
3e^{-\frac{2\pi}{\gamma}} \le 0.0056.
\end{gather*}

In this paper we will extend the above model to its simplest reinforced
random walk form depending on a strength parameter $\de$ and obtain an 
explicit expression for the corresponding $\tP_\de \lf(\la, u,a\rt)$.
Recovering $P_N (x,a|\de)$ from it leads to a very complex result which
we will not bother to write down.  But the same asymptotic limit is a 
different story.  In fact, we observe the equivalence of
$$
N\to\infty \quad\text{at fixed }\quad \frac{a^2}{N} = \gamma \frac{\pi}{2}
+ O\lf(\frac{1}{a}\rt)
$$
and
$$
N= \lf[t/ \ep^2\rt] \qquad a= \lf[b/\ep\rt] \qquad\text{at fixed $t$ 
and $b$}
$$
as $\ep \to 0$ the second of which is precisely the diffusion scaling of
the Brownian motion limit of the walk.  This observation suggest that the
corresponding reinforced Brownian motion may be amenable to detailed 
analysis, and this is what we will do.  It turns out to be simpler
to start from the Brownian Motion version of dynamics than to proceed
via the $g,f$ route, although they must be equivalent  
(and the g.f route has to face the lack of right continuity at the
origin for the limiting Brownian Motion under our left boundary 
condition, and it is well to avoid this).
At any rate, 
we will first define the precise form under consideration and set up
a suitable $g,f$.  We then shift gears and derive the corresponding 
continuous walk.  This is in general represented by an Inverse Laplace
Transform, which attains a simple form in several special cases that 
we obtain.  Finally, we revisit the general $\de >0$ situation and see
what general conclusion can be drawn.

\section{Basic Edge Reinforced Random Walk}

We have already alluded to the tactic of obtaining the statistics
of the maximum of $S_N$ for ordinary next neighbor RW by marking
on the lattice those edges that have been previously visited.  To 
reinforce these edges we now replace them by bold face edges with an 
associated weight $1+\de$ with $\de >-1$, whereas the unreinforced
edges have weight 1.  This recipe is usually called one time 
reinforced \cite{bib5}.  At a given vertex $S_N$, the odds ratio for
the transitions $S_N\to S_N +1$ $S_N\to S_N-1$ is  to be taken as
the ratio of these weights, the possible pairs of transition probabilities
then being $\lf(\half, \half\rt), \lf(\frac{1+\de}{2+ \de},
\frac{1}{2+\de}\rt)$ and $\lf(\frac{1}{2+\de}, \frac{1+\de}{2+\de}\rt)$.
In detail we now study a random walk on the integer lattice 
$x\ge 0$, $a\ge 1$ with joint distribution defined by 
\begin{align*}
P_N (x,a) &= P \{ S_N=x, A_N=a\} \\
\text{where}\quad S_k &= \sum^k_{i=1} X_i\qquad\text{$X_i$ a random variable}
\\[-3mm]
&\hspace{1.25in} X_i = \pm 1
\end{align*}
We deal only with a symmetric random walk reflected at the origin 
(\cite{bib3}, \cite{bib4})
\begin{align*}
P\{ X_i = \pm 1\} = \half \qquad\text{if}\qquad S_{i-1} \neq 0\\
P\{ X_i = 1\} = \half \qquad\text{if}\qquad S_{i-1} = 0\\
\text{and the walk starts at}\quad x=1, a=1
\end{align*}
$$
P_1 (x,a)= \de_{x,1} \, \de_{a,1}.
$$
Note that the initial condition averages equal parity and opposite parity
pairs $(N,x)$, thus avoids the necessity of following the pairs separately.

Since $a$ has not changed from its prior value when $x <a$, we have 
for $x\ge 2$ and $N \ge 1$
\begin{subequations}\label{1}
\begin{align}
P_{N+1} (x,a) &= \half\:P_N (x-1, a) + \half \:P_N (x+1,a)
\quad\text{if}\quad x< a-1
\label{1a}
\\
\begin{split}\label{1b}
P_{N+1} (x,a) &= \half\:P_N (x-2, a) + \frac{1+\de}{2+\de} \:P_N (a,a)
\quad\text{if}\quad x= a-1\\
P_{N+1} (x,a) &= \half\:P_N (x-1, a) + \frac{1}{2+\de} \:P_N (a-1,a-1)
\quad\text{if}\quad x= a\\
P_{N+1} (x,a) &=0 \hspace{2.5in}\text{if} \quad x>a 
\end{split}\\
\begin{split} \label{1c}
P_{N+1} (1,a) &= P_N(0,a) + \half \:P_N (2,a) + \frac{\de}{2(2+\de)}
\:P_N(2,2)\, \de_{a,2}\\
P_{N+1} (0,a) &= \half\:P_N(1,a) + \frac{\de}{2(2+\de)}\:P_N(a,a)\,
\de_{a,1} 
\end{split}
\end{align}
\end{subequations}
Hence for $x\ge 0$\: $a\ge 1$\: $N\ge 1$,
\be\label{2.5}
\begin{aligned}
P_{N+1}(x,a) &= \half \lf(1+ \de_{x,1} - \de_{x,a+1}\rt) P_N(x-1,a)\\
& \;\, + \half\: P_N (x+1,a) + \frac{1}{2+ \de} \:P_N (a-1, a-1)\,
\de_{x,a} (1-\de_{a,1}) \\
&\;\, + \frac{\de}{2(2+\de)} \:P_N (a,a) \, \de_{x,a-1}\\
P_1(x,a) &= \de_{x,1}\; \de_{a,1}.
\end{aligned}
\ee

Our task is now to solve \eqref{2.5}, which we do in standard fashion 
by first introducing the generating function, convergent for $|\la| <1$.
\begin{align*}
\tP (\la, x, a) &= \sum^\infty_{N=1} \la^N  P_N(x,a)\\
&= \la P_1 (x,a) + \sum^\infty_{N=1} \la^{N+1}  P_{N+1} (x,a)
\end{align*}
It follows from \eqref{2.5} that 
\begin{align*}
\tP (\la, x,a) &= \la P_1 (x,a) + \frac{\la}{2} \lf(1+\de_{x,1}
- \de_{x, a+1} \rt) \tP(\la, x-1, a)\\
&\; + \frac{\la}{2} \, \tP (\la, x+1, a) + \frac{\la}{2+ \de}
\:\de_{x,a}\, \tP(\la, a-1, a-1)(1-\de_{a,1})\\
&\; + \frac{\la \de}{2(2+\de)} \: \tP(\la, a, a)\, \de_{x,a-1}.
\end{align*}
Further simplification is then achieved by going over to the double
generating function 
$$
\mP^{\approx} (\la, u,a) \equiv \sum^a_{x=0} \tP(\la, x,a)\,
u^x = \sum^\infty_{N=1} \sum^a_{x=0} \la^N u^x P_N(x,a)
$$
where we have used the fact that $P_N(x,a)=0$ for $x>a$, and this 
also establishes that  $\tP(\la, u,a)$ is a polynomial in $u$ of 
degree $a$, thereby convergent for all $u$. Summing over 
$x$ with weight $u^x$, we find after minor algebra that
\be\label{3}
\begin{aligned}
&\mP^{\approx} (\la, u,a) (2u-\la u^2-\la) -2\la u^2\, \de_{a,1} = \\
&\quad \la (u^2-1) \,\tP(\la,0,a) + \la \,u^a \lf(\frac{\de}{2+\de}-u^2\rt)
\tP(\la, a,a) +\\
&\quad + \frac{2\la}{2+\de} \:u^{a+1} \, \tP (\la, a-1, a-1) (1-\de_{a,1}).
\end{aligned}
\ee
Now let $\tht = \frac{1-\sqrt{1-\la^2}}{\la}$ be the small root of 
$2u-\la u^2-\la=0$; then from \eqref{3} we find
\begin{align}
\lf(1-\tht^2\rt) \la \tP(\la, 0,a) &= 2\la\,\tht^2 \, \de_{a,1}
+ \lf(\frac{\de}{2+\de}-\tht^2\rt) \la \tht^a\,\tP(\la, a, a) \label{4} \\
& \quad + \frac{2\la}{2+\de} \: \tht^{a+1}\,\tP(\la, a-1, a-1) (1-\de_{a,1})
\nonumber\\
\intertext{and}
\lf(1-\frac{1}{\tht^2}\rt) \la \,\tP(\la, 0,a) &= 
\frac{2\la}{\tht^2} \:\de_{a,1} + \frac{\la}{\tht^a} 
\lf(\frac{\de}{2+\de} - \frac{1}{\tht^2}\rt) \tP(\la, a,a) \label{5}\\
&\quad  + \frac{2\la}{2+\de} \:\frac{1}{\tht^{a+1}} \: 
\tP (\la, a-1, a-1) (1-\de_{a,1}).\nonumber
\end{align}
Multiplying \eqref{4} by $\frac{1}{\tht}$ and \eqref{5} by $\tht$
and adding, we get
\begin{align*}
0 &= 2\la(1+\tht^2) \,\de_{a,1} + \lf[ \frac{\de\la}{2+\de}
\lf(\tht^a+ \tht^{2-a}\rt) -\la \lf(\tht^{a+2}-\tht^{-a}\rt)\rt]
\tP (\la, a,a)\\
&\qquad\qquad + \frac{2\la}{2+\de} \lf(\tht^a+\tht^{1-a}\rt) \tP(\la, a-1, a-1)
(1-\de_{a,1})
\end{align*}
or
\begin{align*}
\tP(\la, 1,1) &= \frac{2\tht (1+\tht^2)}{1+\tht^4-
\frac{2\de}{2+\de} \:\tht^2} = \frac{\la}{1-\la^2 \: \frac{1+\de}{2+\de}}
= \sum^\infty_{N=0} \la^{2N+1} \lf(\frac{1+\de}{2+\de}\rt)^N
\\
\tP (\la, a,a) &= \frac{\frac{2}{2+\de} \:\tht\lf(\tht^a+ \tht^{-a}\rt)}{
\tht^{a+2} + \tht^{-a} -\frac{\de}{2+\de} \lf(\tht^a+ \tht^{2-a}\rt)}
\:\tP(\la, a-1,a-1) 
\end{align*}
from which we have
\begin{align}
\tP(\la, a,a) &= \prod^\al_{i=2} \lf[ \frac{1+\tht^{2i}}{1+\tht^{2i+2}
- \frac{\de}{2+\de}\:\tht^2 \lf(1+ \tht^{2i-2}\rt)} \rt] 
\lf(\frac{2\tht}{2+ \de}\rt)^{a-1}\, \tP(\la,1,1)
\label{6} \\
\intertext{and} 
\tP (\la, 0,a) &= \sum^\infty_{N=1} \la^N\,P \lf\{S_N=0, \quad 
\max_{1\le i \le N} \, S_i =a \rt\}\nonumber\\
&= \frac{2\tht^2}{1-\tht^2} \: \de_{a,1} + \frac{\tht^a}{1-\tht^2}
\lf(\frac{\de}{2+ \de}- \tht^2\rt) \tP(\la, a,a) \label{7} \\
&\qquad\qquad + \frac{\tht^{a+1}}{1-\tht^2} \: \frac{2}{2+\de}
\:\tP (\la, a-1, a-1) (1-\de_{a,1}).\nonumber 
\end{align}
Inserting \eqref{6} and \eqref{7} into \eqref{3} we conclude with 
the explicit if involved expression for the double $g, f$:
\be\label{8}
\begin{split}
&\mP^{\approx} (\la, u, a)(2u-\la u^2-\la) - 2\la u^2\, \de_{a,1} =\\
&\quad\la^2 \lf(u^2-1\rt) \frac{2\tht^2}{1-\tht^2}\:\de_{a,1} + 
A\tP (\la, a,a) + B \tP(\la, a-1,a-1)
\end{split}
\ee
where
\begin{align*}
&A = \la \lf(u^2-1\rt)\frac{\tht^a}{1-\tht^a}\lf(\frac{\de}{2+\de}-\tht^2\rt)
+\la u^a \lf(\frac{\de}{2+\de}- u^2\rt)\\
&B = \la\lf(u^2-1\rt) \frac{\tht^{a+1}}{1-\tht^2} \frac{2}{2+\de}
\lf(1-\de_{a,1}\rt) + \frac{2\la}{2+ \de} \: u^{a+1}\lf(1-\de_{a,1}\rt)\\
&\tP (\la, a,a) = \lf(\frac{2\tht}{2+\de}\rt)^{a-1} \tP(\la,1,1)
\prod^a_{i=2} \frac{1+\tht^{2i}}{1+\tht^{2i+2} - \frac{\de \tht^2}{2+\de}
\lf(1-\tht^{2i-2}\rt)}\\
&\quad \text{and }\, \tP(\la,1,1) = \frac{\la}{1-\la^2 \, \frac{1+\de}{2+\de}
}
\end{align*}

\section{Edge Reinforced Brownian Motion}

As previously mentioned, recovering $P_N(x,a)$ from \eqref{8} can 
indeed be carried out, resulting in a very involved multiple summation
that offers no insight into the nature of the process.  However, the
related reinforced Brownian Motion is both solvable, simple in 
special cases, and qualitatively transparent.  We must first convert
the discrete dynamics \eqref{1} into its low resolution continuous
Brownian Motion counterpart.  This requires a rescaling that reduces the
spatial and temporal step size to zero, a highly non-unique procedure.
In Brownian Motion, the 2 rescalings are related in traditional fashion
$$
y= \ep x, \quad b=\ep a, \quad t= \ep^2 N \qquad\text{as}\quad \ep \to 0
$$
resulting formally (see \cite{bib3}) in a 2-dimensional probability density
$$
P(t,y,b) =\lim_{\ep \to 0} \ep^{-2} \,P_{[t/\ep^2]} 
\lf( \lf[\frac{y}{\ep}\rt], \lf[ \frac{b}{\ep}\rt] \rt)
$$
normalized so that
\be\label{9}
\int^\infty_0\int^\infty_0 P(t,y,b) \,dy\, db=1 \qquad\text{for all }\; t>0.
\ee
In what follows, the integer value notation will be implicit, as needed.

We now translate the dynamics given by \eqref{1} into the required 
differential form.  For this purpose, we have already separated \eqref{1}
into 3 regimes, the bulk \eqref{1a}, the right boundary \eqref{1b}, and
the left boundary \eqref{1c}.  We then translate each equation for 
$P_N(x,a)$ into a corresponding equation for $P(t, y,b)$.  First,
\eqref{1a}:
$$
P(t + \ep^2,y, b) - P(t,y,b)= \half \lf\{P(t,y-\ep, b)-2P (t,y,b)+
P(t,y+\ep, b)\rt\}
$$
or dividing by $\ep^2$ and applying l'Hopital,
\begin{subequations}\label{10}
\be
\frac{\p}{\p t}\: P (t,y,b) = \half \:\frac{\p^2}{\p y^2}
\:P(t,y,b) \quad\text{for}\quad y< b.
\label{10a}
\ee
Next, \eqref{1b} which we combine and augment to read
\begin{align*}
&P(t + \ep^2, b-\ep, b) + P(t + \ep^2, b,b) -2P(t,b,b) = \\
&\quad \half \lf[ P(t,b-2\ep, b) + P(t,b-\ep,b) -2P(t,b,b)\rt] +\\
&\quad \frac{1}{2+ \de} \lf[P(t,b-\ep, b-\ep) - P(t,b,b) \rt] 
\end{align*}
Now dividing by $\ep$ and applying l'Hopital as $\ep \to 0$ result in 
\begin{align*}
-\frac{\p}{\p y} P(t,y,b) = - \frac{3}{2} \: \frac{\p}{\p y}
\:P(t,y,b) &- \frac{1}{2+ \de} \lf[  -\frac{\p P(t,y,b)}{\p y} 
- \frac{\p P}{\p b} (t,y,b)\rt]  \\
& \qquad\qquad \qquad
\text{as }\quad y\to b,
\end{align*}
or
\be\label{10b}
\begin{aligned}
&\lf(1+ \frac{\de}{4}\rt) \mP'_y (t,y,b) + \half \, \mP'_b (t,y,b)=0
\quad\text{for}\quad y\to b, \quad\text{and } \\
&\;P(t,y,b)=0 \quad\text{ for }\quad y>b.
\end{aligned}
\ee
Finally \eqref{1c}, where we assume $a>2$ so that $b>0$.  Again, combining
and augmenting, we have 
\begin{align*}
&P(t + \ep^2, \ep, b) +P(t + \ep^2, 0, b) -2P(t,0,b) = \\
&\; \half P (t,2\ep, b) +\half P(t,\ep, b)-P(t, 0, b).
\end{align*}
This approach appears to involve loss of boundary information,
but in fact what it does is to make innocuous any distinction
between differing definitions of reflection at the origin.
At any rate, dividing by $\ep$ and let $\ep\to 0$ now yields, 
on applying l'Hopital,
the expected simple
\be\label{10c}
\mP'_y (t,y,b)=0 \quad\text{for}\quad y\to 0.
\ee
\end{subequations}

To these, we must add the initial condition.  This is most usefully 
regarded as any condition that confines the probability to a finite
$(x,a)$ region on the non-negative half line, since all of these
will under the diffusion scaling reduce to 
\be\label{11}
P(t,y,b) =0 \quad\text{as}\quad t\to 0,
\ee
together with a unit mass at the origin.  The mass is no longer confined
when $t>0$, but replaced by the normalization \eqref{9}.

We must now solve equation \eqref{10a}, \eqref{10b}, \eqref{10c},
\eqref{11} and \eqref{9}.  To do so, it is simplest to work in Laplace
transform
$$
\tP (s,y,b) \equiv \int^\infty_0 e^{-st}\, P(t,y,b)\, dt
$$
and Laplace transform the 3 equations to read
\begin{subequations}\label{12}
\be\label{12a}
s \tP(s,y,b) - \half \,\moP''_ y (s,y,b) =0
\ee
\begin{align}
\lf(1+ \frac{\de}{4}\rt)&\moP'_y(s,y,b) + \half \moP'_b(s,y,b)=0
\quad\text{at}\quad y=b, \;\text{ with}\label{12b}\\
& \tP (s,y,b)=0 \quad\text{ for }\quad y>b \nonumber\\
&\moP'_y (s,y,b)=0 \quad\text{ at} \quad y=0, \label{12c}
\end{align}
\end{subequations}
plus normalization \eqref{9}.

The solution of \eqref{12a} subject to \eqref{12c} is clearly given 
by
$$
\tP (s,y,b) = 
\begin{cases}
f(s,b) \cosh \lf(y\sqrt{2s}\rt) &\text{for  $0 \le y\le b$}\\
0 &\text{for $y>b$}
\end{cases}
$$
and substituting into \eqref{12b}, we see that
$$
\half \frac{f'_b(s,b)}{f(s,b)} + \sqrt{2s} \lf(1+ \frac{\de}{4}\rt)
\tanh \, b\sqrt{2s} =0,
$$
so that
$$
f(s,b)= \frac{k(s)}{\lf(\cosh b\sqrt{2s}\rt)^{2\lf(1+ \frac{\de}{4}\rt)}}\ ,
$$
normalization then evaluating $k(s)$ as $2+ \de$.

Gathering the preceding together we conclude that
\be\label{13}
\tP(s,y,b)= (2+\de) \, 
\frac{\cosh y \sqrt{2s}}{ \lf( \cosh b \sqrt{2s}\rt)^{2 + \frac{\de}{2}} }
\ee

\section{The Marginal Distribution of the maximum}

$P(t,y,b)$ contains a great deal of information, too much perhaps for 
ease of characterization.  A  common first step to reduce the information
overload is to consider separately the state of the Brownian walker, as 
given by 
$$
P(t,y) = \int^\infty_y P(t,y,b)\,db
$$
and that of the lattice
$$
Q(t,b) = \int^b_0 P(t,y,b)\,dy.
$$
For the prototypical case $\de=0$, $P(t,y)$ is that of a primitive
Markovian Brownian Motion, a half-Gaussian, whereas $Q(t,b)$ is more 
involved.  However, it extends very easily to the $\de \neq 0$ case
as follows:

It is sufficient to work directly with the Laplace Transform in time and
then extract the time dependence.  From \eqref{13} we have
\begin{align}
\tQ (s,b) &\equiv \int^b_0 Q(s,y,b) \,dy = \frac{2+\de}{\sqrt{2s}}
\frac{\sinh b\sqrt{2s}}{\lf(\cosh b\sqrt{2s}\rt)^{2+\frac{\de}{2}}} 
\label{14}\\
&= -\frac{2}{b} \frac{\p}{\p s} \lf(\cosh b \sqrt{2s}\rt)^{
-\lf(1+\frac{\de}{2} \rt)} \ .\nonumber 
\\
\intertext{Therefore}
Q(t,b) &= \frac{2t}{b} \mL^{-1}_{t,s} \lf(\cosh b\sqrt{2s}\rt)^{-\lf(1+
\frac{\de}{2}\rt)} \ .
\nonumber
\end{align}
It is ``only'' the technical issue of carrying out the inverse Laplace
Transform $\CL^{-1}_{t,s}$ that remains.  Even for $\de=0$, this requires
introduction of the Jacobi Theta function, not the best way to visualize
the resulting situation.  There are alternatives.  One alternative is to
make use of the rapid convergence of the series expansion
\be\label{15}
\begin{aligned}
\lf(\cosh b\sqrt{2s}\rt)^{-\lf(1+\frac{\de}{2}\rt)} &= 
\lf[\half\lf(e^{b\sqrt{2s}} + e^{-b\sqrt{2s}}\rt) \rt]^{-\lf(1+ 
\frac{\de}{2}\rt)}\\
&= 2^{1+\frac{\de}{2}} \, e^{-\lf(1+\frac{\de}{2} \rt) b\sqrt{2s}}
\lf(1+ e^{-2b \sqrt{2s}}\rt)^{-\lf(1+\frac{\de}{2}\rt)}\\
&= 2^{1+\frac{\de}{2}} \sum^\infty_{j=0} (-1)^j 
\lf(
\begin{matrix}
j+\frac{\de}{2} \\
j
\end{matrix}
\rt)
e^{-\lf(2j+1 + \frac{\de}{2}\rt) b\sqrt{2s}}.
\end{aligned}
\ee
Using $\mL^{-1}_{t,s} \,e^{-c\sqrt{s}} = ce^{-c^2/4t}/2\sqrt{\pi}\,t^{3/2}$,
we therefore have,
\be\label{16}
b\,Q(t,b) =2^{1+\frac{\de}{2}} \, \gamma^{1/2} \sum^\infty_{j=0}(-1)^j
\lf(
\begin{matrix}
j+ \frac{\de}{2}\\
j
\end{matrix} \rt)\!\!
\lf(2j+1+ \frac{\de}{2}\rt)
e^{-\frac{\pi\gamma}{4} \lf(2j+1+\frac{\de}{2}\rt)^2}
\ee
which for sufficiently large $\gamma$ is a rapidly convergent series,
with leading $j=0$ term a pure half Gaussian in $\gamma$, of variance
$\frac{2}{\pi}/1+\frac{\de}{2}$.

If $\gamma$ is small compared to  $\lf(1+\frac{\de}{2}\rt)^{-1/2}$,
the terms in the series increase in absolute value until very large $j$,  
and then start decreasing.  Since convergence is so slow, we would
want to convert to a more appropriate series representation, in the
fashion of the two series referred to in Section~1.  This is carried
out in Appendix A and becomes increasingly complex as $\de$ increases.
Arguably, the  most important information concerning the distribution 
$Q(t,b)$ is that of its moments and cumulants $(\text{coef of } 
\frac{t^k}{k!} \text{ in } \log E(e^{tx}))$.  Let us first obtain 
concise expressions for the moments of $Q(t,b)$.  From \eqref{14}
we have
$$
Q(t,b) = \mL^{-1}_{t,s} \frac{2+\de}{\sqrt{2s}} \frac{\sinh b \sqrt{2s}}{
\lf(\cosh b \sqrt{2s}\rt)^{2+\frac{\de}{2}}}
$$
so assuming the validity of interchanging limits
\begin{align*}
E(b^k) &= \int^\infty_0 b^k \lf( \mL^{-1}_{t,s}
\frac{2+\de}{\sqrt{2s}} \:
\frac{\sinh b \sqrt{2s}}{
\lf(\cosh b \sqrt{2s}\rt)^{2+ \frac{\de}{2}}} \rt)db\\[1mm]
&= \lf(2+\de\rt) \mL^{-1}_{t,s}
\frac{1}{\sqrt{2s}} \int^\infty_0 b^k 
\frac{\sinh b \sqrt{2s}}{\lf(\cosh \sqrt{2s}\rt)^{2+ \frac{\de}{2}}}
\,db
\end{align*}
For $x=b \sqrt{2s}$ this transforms to 
$$
E(b^k)= (2+\de) \,\mL^{-1}_{t,s} \,\frac{1}{(2s)^{\frac{k+2}{2}}}
\int^\infty_0 x^k \frac{\sinh x}{(\cosh x)^{2+ \frac{\de}{2}}}\,dx
$$
Note that $\mL^{-1}_{t,s} \,\frac{1}{s^{p+1}} = \frac{1}{p!} \, t$:
hence
$$
E(b^k) = (2+ \de)\, \frac{1}{2^{\frac{k+2}{2}}} \:t^{\frac{k}{2}}
\frac{1}{\lf(\frac{k}{2}\rt)!} \int^\infty_0 x^k 
\frac{\sinh x}{(\cosh x)^{2+\frac{\de}{2}}} \,dx
$$
Let $y=\cosh x$, we get
$$
E(b^k) = \frac{2+\de}{2} \: \frac{t^{\frac{k}{2}}}{2^{\frac{k}{2}}}
\frac{1}{\lf(\frac{k}{2}\rt)!} \int^\infty_1
\frac{\lf(\cosh^{-1} \,y\rt)^k}{y^{2+\frac{\de}{2}}}
\, dy
$$
For $y=1+z$ and recalling the notation $\de=4m$
\be\label{17}
E(b^k) = (1+2m) \frac{t^{\frac{k}{2}} }{2^{\frac{k}{2}}
\lf( \frac{k}{2}\rt)! } 
\int^\infty_0 
\frac{
\lf[ \cosh^{-1} (1+z)\rt]^k}{(1+z)^{2+2m} }\,
dz
\ee
but for large $m$ it suffices to use the MacLaurin expansion of 
$\cosh^{-1} (1+z)$
$$
\cosh^{-1} (1+z) = \sqrt{2z} \lf(1-\frac{1}{12}z + \frac{3}{40}z^2 + 
\dots \rt)
$$
so that stopping after the first correction term
$$
\lf[ \cosh^{-1} (1+z)\rt]^k = (2z)^{\frac{k}{2}}
\lf(1- \frac{k}{12}z + \dots \rt)
$$
To evaluate \eqref{17} we require the known result
$$
I_k \equiv \int^\infty_0 \frac{z^{k/2}}{(1+z)^{2m+2}}
= \frac{
\lf(2m-\frac{k}{2}\rt)! \lf(\frac{k}{2}\rt)!}{
\lf(2m+1\rt)!}
$$
Hence
\begin{align*}
E(b^k) &= 
(1+2m) \,t^{\frac{k}{2}} \frac{I_k}{\lf(\frac{k}{2}\rt)!}
\lf[ 1-\frac{k}{24} \,\frac{I_{k+2}}{I_k}\rt] \\
&= (1+2m)\, t^{\frac{k}{2}} \, \frac{\lf(2m-\frac{k}{2}\rt)!}{
\lf(2m+ 1\rt)!}
\lf[ 1-\frac{k}{24} \,\frac{\frac{k}{2} +1}{2m-
\frac{k}{2}} \rt]
\end{align*}
One can show after some algebra that for large $m$
$$
\frac{\lf(2m -\frac{k}{2}\rt)!}{\lf(2m+1\rt)!} =
\frac{1}{\lf(2m+1 \rt)^{\frac{k}{2} +1}}
\lf( 1+ \frac{1}{2m+1} \,\frac{k}{4} \lf( \frac{k}{2} +1 \rt)
+ \dots \rt)
$$
It follows that
$$
E(b^k) = \frac{t^{k/2}}{\lf(2m+1\rt)^{k/2}} 
\lf( 1 + \frac{5k}{24} \: \frac{\lf(\frac{k}{2} +1\rt)}{
\lf(2m+1\rt)} + \dots \rt)
$$
In particular to leading order in $(2m+1)$
$$
E(b) = \frac{t^{1/2}}{\lf(2m+1\rt)^{1/2}}
$$
whereas
$$
\sigma^2_b = \text{Var}\,b = E(b^2) - \lf(E(b)\rt)^2 = \frac{5}{24}
\:\frac{t}{\lf(2m+1\rt)^2}
$$
so that
$$
\frac{\sigma_b}{E(b)} = \sqrt{\frac{5}{24(2m+1)}}.
$$

\section{The Marginal Distribution of the Walker}

The general behavior of $Q_m(t,b)$ and its moments depends quantitatively
on the value of $\de \ge 0$ but not qualitatively (Theta related functions
are always encountered).  On the other hand, $P_m(t, y)$ the walker
distribution, is particularly simple in form when $\de =0$, $(m=0)$ but 
one should not expect this simplicity to be maintained for $\de >0\, (m>0)$.
Let us see how its form and its moments are altered.

From \eqref{13}
\be\label{18}
\moP_m (s,y,b) = (2+\de) \,\frac{\cosh y \sqrt{2s}}{
\lf(\cosh b \sqrt{2s}\rt)^{2+\frac{\de}{2}}}\ ,  \qquad \de=4m
\ee
From \eqref{18} 
$$
P(t,0)=\CL^{-1}_{t,s} \lf(1+2m\rt)
2\int^\infty_0
\frac{db}{(\cosh \sqrt{2s})^{2m+2}} = \sqrt{\frac{2}{t\pi}}
\,\frac{2^{2m}}{(2m)!} \, (m!)^2
$$
and asymptotically for large $m$ $P(t,0)=\sqrt{\frac{2m+1}{t}}$.

Furthermore from \eqref{15}
\begin{align*}
&(2+\de) \,\cosh \,y\sqrt{2s} \lf(\cosh b \sqrt{2s}\rt)^{
\lf(2+ \frac{\de}{2}\rt)} = \\
&\qquad\qquad 2^{2+ \frac{\de}{2}} \,(2+\de)
\sum_{j=0} (-1)^j \lf( \begin{matrix}
j+1 + \frac{\de}{2} \\
j 
\end{matrix} \rt)\!
e^{-\lf(2j+2+ \frac{\de}{2}\rt)b\sqrt{2s}}
\; \cosh y \sqrt{2s} 
\end{align*}
It follows that
\begin{align*}
&\moP_m (s,y) = \int^\infty_y \moP_m (s,y,b)\,db = \\
& (2+\de)\,2^{1+\frac{\de}{2}}\sum^\infty_{j=0} (-1)^j
\lf(\begin{matrix}
j+1 + \frac{\de}{2}\\
j \end{matrix} \rt)
\frac{e^{-y\sqrt{2s}\lf(1+2j+ \frac{\de}{2}\rt)} + 
e^{-y \sqrt{2s} \lf(3+2j+ \frac{\de}{2}\rt)} }{
\sqrt{2s} \lf(2j +2 + \frac{\de}{2}\rt)} \nonumber
\end{align*}
Carrying out the inverse Laplace transform we get 
\begin{align*}
\mP_m(t,y) &= \frac{\lf(2+4m\rt)2^{2m+1}}{\sqrt{2\pi \,t}}
\sum_{j=0} (-1)^j \lf(\begin{matrix}
j+1 + 2m\\
j\end{matrix} \rt) 
\lf[ e^{-\frac{y^2}{2t}\lf(2j+2m+1\rt)^2} + \rt. 
\\
&\qquad\qquad \lf. e^{\frac{-y^2}{2t} \lf(2j+2m+3\rt)^2} \rt] 
\frac{1}{2j+2+2m} 
\end{align*}

Since the two series converge absolutely we can shift the index of one
of them and readily obtain:
\be\label{19}
\begin{aligned}
&\mP_m (t,y)  = \frac{2+2m}{\sqrt{2\pi t}} \, 2^{1+2m}\,
\lf\{ 
\frac{1}{2+2m} \,e^{-\frac{y2}{2t}\lf(1+2m\rt)^2} -
\sum^\infty_{j=0} (-1)^j \,C_j \,e^{-\frac{y^2}{2t} \lf(3 + 2j + 2m\rt)^2}
\rt\} \\
&\qquad\qquad\text{where }\quad  C_j = \frac{2m}{(1+2m)}
\frac{(2j+3+2m)}{\lf[(2j+3+2m)^2-1\rt]}
\lf(\begin{matrix}
j+1 + 2m\\ 
2m \end{matrix} \rt)
\end{aligned}
\ee
Equation \eqref{19} is an alternating series, and for $\frac{y^2}{t}
> \frac{\log 2(m+1)}{4(m+1)}$ it is monotonically decreasing in magnitude.
\eqref{19} is therefore dominated by its first term, which is the exact
result when $m=0$.  Therefore, when 
\be\label{20}
\frac{y^2}{t} > \frac{\log 2(m+1)}{4(m+1)}
\ee
there is no qualitative change for $m>0$.

When \eqref{20} is not satisfied, convergence is slow, but in principle
we can convert  the series in \eqref{19} into one in which $\frac{t}{y^2}$
appears in the exponent.  The result is very complicated.  Therefore, 
instead of finding the full distribution we will only find $E(y)$ and
$E(y^2)$.  We will do this by direct integration of the series in \eqref{19}.
To allow the interchange of the infinite sum with the improper integral,
we will first insert a convergence factor $\al^j$ into \eqref{19}, later 
to be set to $\al=1$.

\subsection*{a.  The Mean}
\begin{align*}
\text{For  $m\neq 0$ } \quad E(y) &= \int^\infty_0 y \mP_m (t,y)\,dy\\
&= \frac{2+4m}{\sqrt{2\pi t}} 2^{1+2m} 
\lf\{\rule{0in}{.27in}
\frac{1}{2(m+1)} \int^\infty_0 y e^{-\frac{y^2}{2t}(1+2m)^2}\, dy \rt. \\
&\qquad\qquad\qquad  \lf. - \sum_{j=0} (-1)^j \,C_j \al^{j+1} \int^\infty_0
y e^{-\frac{y^2}{2t}(3+2j+2m)^2}\, dy \rule{0in}{.27in}\rt\} \\
&= \frac{(1+2m)\,2^{2+2m}}{\sqrt{2\pi t}} \lf\{\rule{0in}{.27in}
\frac{1}{2(m+1)} \frac{t}{(1+2m)^2} \rt.  \\
&\qquad\qquad\qquad  \lf. - \sum_{j=0} (-1)^j \,C_j \frac{t\,\al^{j+1}}{
(3+2j+2m)^2}\rule{0in}{.27in}\rt\} \\
&= \frac{2^{2(m+1)}\,t}{\sqrt{2\pi t}} \sum^\infty_{j=0} (-1)^j\,\al^j\,
\frac{\Gamma  (2m+2j) \Gamma(3)}{\Gamma(2m+2j +3)} 
\; m\left(\begin{matrix}
j+2m\\ 2m \end{matrix} \rt)
\end{align*}
Hence (valid for $m=0$ as well, since $lim_{m\to 0} B(2m,3) m\to 1$)
\begin{align*} 
E(y)_\al &= \frac{2^{2(m+1)}\,t}{\sqrt{2\pi t}}
\sum^\infty_{j=0} (-1)^j\,\al^j\, B
\lf(2m+2j,3\rt) 
m\lf(\begin{matrix}
j+2m\\ 2m\end{matrix}\rt) \\
&=\frac{2^{2(m+1)}\,t}{\sqrt{2\pi t}}
\sum^\infty_{j=0} \int^1_0 (-1)^j \, \al^j \,
m\lf(\begin{matrix}
j+2m\\ 2m \end{matrix} \rt)
u^{2m+2j-1}\lf(1-u\rt)^2 du\\
&=\frac{2^{2(m+1)}\,t}{\sqrt{2\pi t}}
m \int^1_0 \lf(1-u\rt)^2 u^{2m-1} \lf(1+\al u^2\rt)^{-(2m+1)}\, du
\end{align*}
or
\be\label{21}
\mE_m (y| \al=1) = \frac{2^{2(m+1)}}{\sqrt{2\pi t}}
tm \int^1_0 (1-u)^2\, u^{2m-1} \lf(1+u^2\rt)^{-(2m+1)} \, du
\ee
\eqref{21} is readily evaluated e.g.\ by setting $u= \tan \frac{x}{2}$ so that
\begin{align*}
\mE_m (y) &= m\sqrt{\frac{t}{2\pi}} \int^1_0 \lf(\frac{1-u}{u}\rt)^2
\lf(\frac{2u}{1+u^2}\rt)^{2m+1} \, du \\
&= m\sqrt{\frac{t}{2\pi}}\: 2 \int^{\frac{\pi}{2}}_0 
\lf[ \lf(\sin x\rt)^{2m-1} - \lf(\sin x\rt)^{2m}\rt]dx \\
&= m \sqrt{\frac{t}{2}} \lf\{ 
\frac{(m-1)!}{\lf(m-\half\rt)!} - \frac{\lf(m-\half\rt)!}{m!}\rt\}.
\end{align*}
Since $\lf(m-\half\rt)!=\frac{(2m)! \,\sqrt{\pi}}{2^{2m}\, m!}$, then 
$$
\text{For asymptotically large} \quad  m,\qquad  
\mE_m (y) \sim \frac{1}{4}
\sqrt{\frac{t}{2m}} \ .
$$
\subsection*{b. Second Moment}

In the same fashion, 
\begin{align*}
\mE_m (y^2) &= \int^\infty_0 y^2 \mP_m (t,y)\,dy\\
&= \frac{2+4m}{\sqrt{2\pi t}} \; 2^{1+2m}
\lf\{\rule{0in}{.27in} 
\frac{1}{2(m+1)} \int^\infty_0 y^2 \,
e^{-\frac{y^2}{2t} (1+2m)^2 }\, dy \rt. \\
& \qquad\qquad  + \lf.  
\sum_{j=0} (-1)^{j+1} \, \al^{j+1}\, C_j 
\int^\infty_0  y^2\,e^{-\frac{y^2}{2t}\,(3+2j+2m)^2}\, dy 
\rule{0in}{.27in}\rt\}
\end{align*}
or carrying out the integrations and inserting the value of $C_j$ (see 
\eqref{19}) we get
\be\label{22}
\begin{aligned}
\mE_m(y^2) &= (1+2m)\, 2^{2m+1}\, t 
\lf\{ \sum_{j=0} (-1)^j \, \al^j \times 
\lf(\begin{matrix}
j+2m\\ 2m \end{matrix} \rt) \times \rt. \\
& \lf.  \qquad \frac{2m}{(2m+1)(2m+2j) (2m+2j+1)^2 \,(2m+2j+2)} 
\rule{0in}{.27in}\rt\}
\end{aligned}
\ee
The $j$ dependence of the last term can be represented as a product of 
two Beta functions, so that \eqref{22} becomes
\be\label{23}
\begin{aligned}
\mE_m(y^2) &= 
(1+2m)\,2^{2m+1}\, t
\sum^\infty_{j=0} (-1)^j\, \al^j\, \frac{2m}{2m+1} 
\lf(\begin{matrix} 
j+2m\\ 2m \end{matrix} \rt) 
\int^1_0\int^1_0
\lf\{\rule{0in}{.27in} \rt.  \\
&\qquad\qquad\qquad  \lf. x^{2j}\, z^{2j} (1-x) (1-z) x^{2m-1} z^{2m} \,dxdz
\rule{0in}{.27in}\rt\} \\
&=(2m+1)\, 2^{2m+1}\,t \int^1_0 \int^1_0 \lf(1+ x^2\,z^2\rt)^{-(2m+1)}
(1-x) (1-z) \, x^{2m-1} \,z^{2m}\, dxdy.
\end{aligned}
\ee
\eqref{23} cannot, for arbitrary $m$, be expressed in terms of 
elementary functions but we can do so for both small $m$ and large $m$.

For small $m$, we find (using Mathematica): 
\begin{alignat*}{2}
\mE_1(y) &= \frac{2-\frac{\pi}{2}}{\sqrt{2\pi}} = 0.171227 &
\mE_1(y^2) &= 3(1-G)= 0.25210  \\[-2mm]
&&&\text{ where $G$ is Catalan's constant.}\nonumber\\
\mE_2(y) &= \lf(\frac{4}{3} - \frac{3\pi}{8}\rt) \sqrt{\frac{2}{\pi}}
= 0.12386 &\qquad
\mE_2(y^2) &= \frac{5}{119} \lf(68-21\pi\rt)  \\
&&& = 0.0703665
\end{alignat*}
For large $m$ we readily find that
\begin{align*}
&\lim_{m\to\infty} (2m+1) \, \mE_m(y^2) = \\
&\quad \lim_{m\to\infty} (2m+1)^2\, t \int^1_0\int^1_0
\lf(\frac{2xz}{1+x^2z^2}\rt)^{2m+1} (1-x)(1-z) \, \frac{1}{x^2}
\:\frac{1}{z} \, dxdz = \\ 
&\qquad \qquad t \int^\infty_0\int^\infty_0  uv e^{-\half (u+v)^2}
\, dudv=\frac{t}{3}
\end{align*}
It follows that asymptotically,
$$
\mE_m(y^2) \sim \frac{1}{3} \; \frac{t}{2m+1}\ .
$$
Note that this asymptotic relationship remains very accurate down to $m=2$.

The variance of $y$ is
\begin{align*}
&\mE_m(y^2) - \lf[ \mE_m(y)\rt]^2 = \frac{1}{3} \: \frac{t}{2m+1} -
\frac{1}{16}\: \frac{t}{2m+1} \\
&\hspace{1.15in}= \frac{13}{48}\:\frac{t}{2m+1} \\
\text{with}\quad &\frac{\sigma_y}{E(y)} =  \frac{\sqrt{\frac{13}{48}}}{
\frac{1}{4}} = \sqrt{\frac{13}{3}} \sim 2.08
\end{align*}
a consequence of the fact that the walker still spends most of its
time around the origin.

\section{Concluding Remarks}

We have examined in some detail the statistics of reinforced Brownian
Motion on the half line, a generalization of the primitive random walk,
in which each edge traversed is given an enhanced weight of $1+\de$.
We did this by applying the diffusion scaling limit to discrete walker
dynamics, resulting in two dimensional dynamics which is diffusive only 
in one dimension, but coupled by oblique reflection at one boundary.
The principal qualitative change in the walker distribution is a contraction
towards the origin.  This augurs well for the study of the much more 
complex case in which the enhanced weight depends upon the number of times
an edge has been traversed.  We are now attending to this system.

\appendix

\section*{Appendix A}

The analysis simplifies materially if $\de=2n$ for $n$ integer
which we henceworth assume.

In \eqref{16} we make a replacement $j\to -1-n-j$ so that \eqref{16} implies
\be\tag{A1}
b\,Q(t,b) = 2^{1+n}\,\gamma^{1/2} \sum^{-(n+1)}_{j=-\infty} (-1)^j
\lf(
\begin{matrix}
j+n\\
n
\end{matrix} \rt)\!
(2j+1+n)\,e^{-\frac{\pi\gamma}{4} \lf(2j+1+n\rt)^2}
\ee
Note that $\lf(\begin{smallmatrix} j+n\\n\end{smallmatrix}\rt)=0$ whenever
$-n <j<0$ therefore the upper limit in \thetag{A1} can be replaced by $-1$.
Doing so and adding \eqref{16} to \thetag{A1} we get
\be\tag{A2}
2b Q(t,b)=2^{1+n} \,\gamma^{1/2} \sum^\infty_{j=-\infty} (-1)^j
\lf(\begin{matrix}
j+n\\ n \end{matrix} \rt)\!
(2j+1+n) \, e^{-\frac{\pi\gamma}{4}\lf(2j+1+n\rt)^2}.
\ee
We can now apply the extended Poisson Transformation which takes the form
(see \cite{bib2})
\begin{align}
&\sum^\infty_{j=-\infty} (-1)^j \,f(2j+1) = \frac{1}{2i}
\sum^\infty_{j=-\infty}(-1)^j
\,\bF \lf(\frac{1}{4} (2j+1)\rt) \tag{A3}\\
\intertext{where}
&\bF(j) \equiv \int^\infty_{-\infty} f(w)\,e^{2\pi ijw}\,dw.\tag{A4}
\end{align}

In order to apply \thetag{A4} to \thetag{A2} we 
assume $n=2m$ i.e.\ $m=\frac{\de}{4}$
$$
Q_m(t,b) \equiv 2^{2m}\,\gamma^{\half} \sum^\infty_{j=-\infty} (-1)^j
\lf(\begin{matrix}
j+2m\\ 2m \end{matrix}\rt)\!
(2j+1+2m)\, e^{-\frac{\pi \gamma}{4} \lf(2j+1+2m\rt)^2}.
$$
Replace $j$ by $j-m$; then 
\be
b Q_m (t,b) =2^{2m} \,(-1)^m\, \sqrt{\gamma} \sum^\infty_{j=-\infty}
(-1)^j \lf(\begin{matrix}
j+m\\ 2m
\end{matrix} \rt)\!(2j+1)\, e^{-\frac{\pi\gamma}{4}\lf(2j+1\rt)^2}
\tag{A5}
\ee
In particular $b Q(t,b|\de=0) = \sqrt{\gamma} \, A(\gamma)$
\be\tag{A6}
\text{where }\, A(\gamma) \equiv \sum^\infty_{j=-\infty} (-1)^j
(2j+1)\,e^{-\frac{\pi\gamma}{4} \lf(2j+1\rt)^2} = \gamma^{-\half}
\,Q_0 (t,b).
\ee
In the notation of \thetag{A3} $f(w)$ belonging to $A(\gamma)$ is 
given by 
$$
f(w)  = w \,e^{-\frac{\pi\gamma}{4}\, w^2}
$$
Carrying out the Fourier transform \thetag{A3} now implies
\be\tag{A7}
A(\gamma) = \frac{2}{\gamma^{3/2}} \sum^\infty_{j=-\infty}
(2j+1)(-1)^j \, e^{-\frac{\pi}{4\gamma}\lf(2j+1\rt)^2}
\ee
For $\de \neq 0$ we will need
\begin{align*}
A_p(\gamma) &\equiv \sum^\infty_{j=-\infty} (-1)^j (2j+1)^{2p+1}\,
e^{-\frac{\pi\gamma}{4}\lf(2j+1\rt)^2} \\
&= \lf(\frac{4}{\pi}\rt)^P \lf(\frac{\p}{\p \gamma}\rt)^P
(-1)^P\, A(\gamma)
\end{align*}
As an example consider $\de=4$ $(m=1)$.
Then from \thetag{A5}
\be\tag{A8}
\begin{aligned}
b Q_1(t,b) &\equiv bQ\lf(t,b| \de=4\rt) = 4\sqrt{\gamma} 
\sum^\infty_{j=-\infty} (-1)^{j+1}\!\lf(\begin{matrix}
j+1\\ 2 \end{matrix} \rt)\!
(2j+1)\,e^{-\frac{\pi\gamma}{4}\lf(2j+1\rt)^2} \\
&= \sqrt{\gamma}\, \half \sum^\infty_{j=-\infty} (-1)^j
\lf[2j+1- (2j+1)^3\rt] e^{-\frac{\pi \gamma}{4} \lf(2j+1\rt)^2} \\
&= \frac{\sqrt{\gamma}}{2} \lf[ A(\gamma) - A_1(\gamma)\rt]
= \frac{\sqrt{\gamma}}{2} \lf[ A(\gamma) + \frac{4}{\pi}\;
\frac{\p}{\p \gamma}\: A(\gamma)\rt]
\end{aligned}
\ee
Using the equality of \thetag{A6} and \thetag{A7} together with 
\thetag{A8} we get
\be\tag{A9}
\begin{aligned}
b Q_1(t,b|\de=4) &= \frac{1}{\gamma} \sum^\infty_{j=-\infty}
(-1)^j (2j+1)\,e^{-\frac{\pi}{4\gamma}\lf(2j+1\rt)^2}\\
&- \frac{6}{\pi\gamma^2} \sum^\infty_{j=-\infty} (-1)^j
(2j+1)\,e^{-\frac{\pi}{4\gamma} \lf(2j+1\rt)^2} \\
&+ \frac{1}{\gamma^3} \sum^\infty_{j=-\infty} (2j+1)^3(-1)^j\,
e^{-\frac{\pi}{4\gamma} \lf(2j+1\rt)^2}
\end{aligned}
\ee
Extending \thetag{A9} to higher values of $\de$ becomes increasingly 
complex.  


\newpage




\end{document}